\documentclass[10pt]{article}

\usepackage{amsthm}   
\usepackage{amssymb}
\usepackage{amsmath}
\theoremstyle{plain}
\newtheorem{app}{Application}[section]
\newtheorem{corollary}[app]{Corollary}
\newtheorem{lemma}[app]{Lemma}
\newtheorem{proposition}[app]{Proposition}
\newtheorem{theorem}[app]{Theorem}

\theoremstyle{definition}
\newtheorem{definition}[app]{Definition}
\newtheorem{example}[app]{Example}

\newtheorem*{note}{Note}

\theoremstyle{remark}
\newtheorem{remark}{Remark}[section]
\newtheorem*{uremark}{Remark}

\newcommand{\bN}{\mathbb{N}}
\newcommand{\bK}{\mathbb{K}}
\newcommand{\bP}{\mathbb{P}}
\newcommand{\bQ}{\mathbb{Q}}
\newcommand{\bR}{\mathbb{R}}
\newcommand{\bZ}{\mathbb{Z}}

\newcommand{\cB}{\mathcal{B}}

\newcommand{\cK}{\mathcal{K}}
\newcommand{\cL}{\mathcal{L}}
\newcommand{\cM}{\mathcal{M}}
\newcommand{\cN}{\mathcal{N}}

\newcommand{\Qsym}{\text{QSym}}
\newcommand{\Mat}{\text{Mat}}
\newcommand{\frakm}{\mathfrak{m}}
\newcommand{\xx}{\mathbf{x}}

\newcommand{\refines}{\preccurlyeq}   
\newcommand{\conv}{\text{conv}}
\newcommand{\spam}{\text{span}}
\newcommand{\supp}{\text{supp}}
\newcommand{\fatzero}{\mathbf{0}}
\newcommand{\disjointsum}{+}    
\newcommand{\suchthat}{:}

\newcommand{\permsof}{S}  

\newcommand{\BJR}{Billera, Jia, and Reiner }

\begin{document}

\title{A Matroid-Friendly Basis for the Quasisymmetric Functions}
\author{Kurt W. Luoto \\
Department of Mathematics, University of Washington, \\ 
Box 354350, Seattle, WA 98195, USA \\
kwluoto@math.washington.edu
}
\maketitle

\begin{abstract}
A new $\bZ$-basis for the space of quasisymmetric functions ($\Qsym$, for short) is presented.
It is shown to have nonnegative structure constants, and several interesting properties relative to
the  quasisymmetric functions associated to matroids by the Hopf algebra morphism $F$ of Billera, Jia, and Reiner \cite{billera-jia-reiner}.
In particular, for loopless matroids, this basis reflects the grading by matroid rank, as well as by the size of the ground set.
It is shown that the morphism $F$ distinguishes isomorphism classes of rank two matroids, and that decomposability of the quasisymmetric function of a rank two matroid mirrors the decomposability of its base polytope.
An affirmative answer to the Hilbert basis question raised in \cite{billera-jia-reiner} is given.
\end{abstract}

\section{Introduction}

In this paper we construct a new $\bZ$-basis for the space of quasisymmetric functions, $\Qsym$
and study its properties. 
For instance, we show that it has nonnegative structure constants, 
and that it behaves well with respect to 
the quasisymmetric functions associated to matroids by the Hopf algebra morphism $\Mat \to \Qsym$ described by \BJR \cite{billera-jia-reiner}.
We also answer in the affirmative a question regarding rank two matroids posed in  \cite[Question 7.10]{billera-jia-reiner}, and give an affirmative answer to  \cite[Question 7.12]{billera-jia-reiner} in the case of rank two matroids.

In  \cite{billera-jia-reiner}, \BJR describe an invariant for matroids in the form of a quasisymmetric function.
They show that the mapping  $F: \Mat \to \Qsym$ is in fact a morphism of combinatorial Hopf algebras (given a suitable choice of character on $\Mat$; see \cite{aguiar-bergeron-sottile}), 
where $\Mat$ is the Hopf algebra of matroids introduced by Schmitt \cite{schmitt-I},  and studied by Crapo and Schmitt 
\cite{crapo-schmitt-4}, \cite{crapo-schmitt-1}, \cite{crapo-schmitt-2}, \cite{crapo-schmitt-3}.
\BJR  show that, while the mapping $F$ is not surjective over integer coefficients, it is surjective over rational coefficients.

Our new basis for $\Qsym$ is ``matroid-friendly" in that it reflects the rank of loopless matroids as well as the size of the ground sets:
for every $1 \leq r \leq n$, there is a set $\cN^n_r$ of ${n-1 \choose r-1}$ basis vectors such that for every loopless matroid $M$ of rank $r$ on an $n$-element ground set, $F(M) \in \spam\; \cN^n_r$;
moreover, $\Qsym$ decomposes as the direct sum of these subspaces.
This provides us with a new product grading of $\Qsym$, according to matroid rank $r$.
(The usual grading of $\Qsym$ by degree corresponds to the size $n$ of the matroid ground set.)
Also, as with the monomial and fundamental bases of $\Qsym$, for every matroid $M$, $F(M)$ has nonnegative coefficients in our basis.


\bigskip

The paper has two main parts.
The first part (Sections \ref{sec:prelims}--\ref{sec:F-of-p}) presents the new basis and relevant background material.
In Section \ref{sec:prelims},
we recount background material  from the literature regarding posets and quasisymmetric functions.
In Section \ref{sec:new-basis}, 
we present a definition for our new basis for $\Qsym$ by means of a construction, and highlight several of its important features.
There we also prove that it is a $\bZ$-basis for $\Qsym$.
In Section \ref{sec:F-of-p}, 
we build necessary machinery regarding computing the quasisymmetric function associated to a labeled poset, in the form of alternative decompositions,
and apply these tools to prove that the structure constants of the new basis are nonnegative.

The second part, (Sections \ref{sec:matroids}--\ref{sec:observations})
discusses matroids and their quasisymmetric functions.
In Section \ref{sec:matroids}, 
we recall some of the concepts, terminology, and results from the paper  \cite{billera-jia-reiner},
and prove our claims regarding the quasisymmetric functions of matroids vis-a-vis our new basis.
In Section \ref{sec:rank2},
we recall  the context of \cite[Section 7] {billera-jia-reiner} regarding
the relationship between decompositions of the quasisymmetric function associated to a matroid
and decompositions of its matroid base polytope,
and recall the statement of \cite[Question 7.10] {billera-jia-reiner} regarding the functions associated to rank two matroids.
We develop a formula for the quasisymmetric function of a loopless rank two matroid in terms of the new basis, and apply it to show
(1) that the morphism $F: \Mat \to \Qsym$ distinguishes isomorphism classes of rank two matroids,
(2) that the two types of decompositions mirror each other, i.e. an affirmative answer to \cite[Question 7.12] {billera-jia-reiner}  for the case of rank two matroids,
and (3) to give an affirmative answer to  \cite[Question 7.10] {billera-jia-reiner}.
In Section \ref{sec:observations}, 
we make additional observations regarding matroid functions and the new basis.
We also compare the new basis with the other $\Qsym$ bases discussed in Section 10 of \cite{billera-jia-reiner}, and sketch an alternate proof of the surjectivity of the map $\Mat \to \Qsym$ over rational coefficients.

\section{Preliminaries}\label{sec:prelims}

In this section we quote certain concepts, terminology, and facts from the literature,
as well as establish certain conventions which will be used in the remainder of the paper.
 
\subsection{Compositions}

A \emph{composition} $\alpha$ is a finite sequence of positive integers,
i.e. $\alpha \in \bP^m$  for some $m \in \bN$.
The number of \emph{parts} of $\alpha$, $m$,
is the \emph{length} of $\alpha$, and denoted by $\ell(\alpha)$.
The \emph{weight} of $\alpha = (\alpha_1,\ldots,\alpha_m)$ is $|\alpha| = \alpha_1+\cdots+\alpha_m$.
Included in our definition is the composition having no parts, 
which we denote by the (bold font) symbol $\fatzero$.
We have $\ell(\fatzero) = |\fatzero| = 0$, the only composition with these properties.

Note, for small examples where individual parts are less than 10, we will often write a composition as a sequence of digits, with no separating commas.
For example, we may write $(1,5,6,3,2,3)$ as $156323$ when the context is clear.
We adopt a similar convention for the one-line notation of permutations in $S_n$ when $n < 10$.

There is a natural bijection between compositions  of weight $|\alpha| = n$ and susbets of $[n-1]$
(where $[n] = \{1,2,3,\ldots,n\}$), 
given by
\[ (\alpha_1,\ldots,\alpha_m) \leftrightarrow \{\alpha_1,\alpha_1+\alpha_2,\alpha_1+\alpha_2+\alpha_3,\ldots,\alpha_1+\cdots+\alpha_{m-1}\}. \]

We say that $\beta$ is a \emph{refinement} of $\alpha$, or that $\beta$ refines $\alpha$ (denoted $\beta \refines \alpha$) if $|\alpha| = |\beta|$ and $A \subset B$,
where $A$ and $B$ are the sets associated to $\alpha$ and $\beta$ respectively.

\bigskip
To any permutation $\pi \in S_n$ there is an associated composition of weight $n$ which we denote $C(\pi)$ and whose parts give the lengths of successive increasing runs in the one-line notation of $\pi$.
For example, for $\pi = 934756218 \in S_9$, we have $C(\pi) = 13212$. 
In this paper, we mildly generalize the notion of a permutation to be any sequence of distinct positive integers.
Given a set of positive integers $X$, we let $\permsof(X)$ denote the set of all permuations of all the elements of $X$.
The run length operator $C(\pi)$ extends to these general permutations in the obvious way.
If $X$ and $Y$ are two sets of positive integers of the same cardinality $n$,
then every bijection $f:X \to Y$ induces a mapping $f:\permsof(X) \to \permsof(Y)$ given by $f(x_1,\ldots,x_n) = (f(x_1),\ldots,f(x_n))$.
If $f$ is an increasing function, then we have $C(f(\pi)) = C(\pi)$ for every $\pi \in \permsof(X)$.

\subsection{Well-known $\Qsym$ bases}

The algebra of quasisymmetric functions $\Qsym$
(or $\Qsym(\xx)$ when we want to emphasize the variable set)
forms a subring of the power series ring $R[[\xx]]$ 
where $\xx = (x_1, x_2, x_3, ... )$ is a linearly ordered set of variables indexed by the positive integers,
and $R$ is a (fixed) commutative ring.
In this paper we only deal with the cases where $R$ is either $\bZ$ or $\bQ$, assuming coefficients in $\bQ$ unless otherwise stated.
We often suppress the variables in our notation, 
writing simply $f \in \Qsym$ rather than $f(\xx) \in\Qsym(\xx)$.

There are a number of well-known bases for $\Qsym$, all indexed by compositions.
(For the two considered here, see \cite{gessel}.)
The best-known is the basis of \emph{monomial} quasisymmetric functions,
which here we denote $\{ \xx^\alpha \}$. 
Given a composition $\alpha$ with $\ell(\alpha) = k$, $\xx^\alpha$ is defined by
\[ \xx^\alpha := \sum_{ 1 \leq i_1 < i_2 < \cdots < i_k } x_{i_1}^{\alpha_1} x_{i_2}^{\alpha_2} \cdots x_{i_k}^{\alpha_k}. \]

Another frequently used basis is the set $\{ L_\alpha \}$ of \emph{fundamental} quasisymmetric functions, defined by \[  L_\alpha := \sum_{\beta \refines \alpha} \xx^\beta. \]
For example, $L_1 = \xx^1$ is simply the degree one elementary symmetric function.

We note that $\Qsym$ as an algebra under the usual multiplication is graded by degree.
For each of the bases described above, the set of basis elements indexed by all the compositions of a fixed weight $n$ forms a basis for the homogeneous component of degree $n$, $\Qsym_n$.
Accordingly, $\dim \Qsym_n = 2^{n-1}$.

\subsection{Posets and $P$-partitions}

One of the early references to quasisymmetric functions is the paper of Gessel \cite{gessel}
(who built on the work of Stanley \cite{stanley-thesis}),
where they are related to $P$-partitions of labeled posets. 
Most of this material can also be found in Stanley \cite{stanley-ec2}.
In the following, we let $\leq$ denote the usual ordering on integers, 
and $\leq_P$ denote the partial order of a poset $P$.
All posets we consider here are finite.

We adopt a mild generalization of Gessel's convention.
We say that a \emph{labeled poset} on $n$ elements is 
a partial order on a set of $n$ positive integers.
These integers are referred to as the \emph{labels} of the poset. 
Usually the set of labels is $[n] = \{1,2,\ldots,n\}$,
but sometimes we make use of other labels.
We often use the same symbol to refer to both the poset and its set of labels when the meaning is clear from context.

\begin{note}
This convention differs from that used in  \cite{billera-jia-reiner}.
There, a labeled poset consists of a pair $(P,\gamma)$,
where $P$ is a poset on an arbitrary set of $n$ elements, 
and $\gamma$ is a \emph{labeling} of $P$,
that is a bijection between the elements of P and the set $[n]$.
The notion is equivalent to Gessel's.
For our generalization, the labeling would be an injective function from the set of elements of the poset into the set $\bP$ of positive integers.
At times we find it convenient to write $(P,\gamma)$
when we wish to discuss various labelings on the same underlying unlabeled poset.
\end{note}

The following is not the actual definition used in \cite{gessel} and \cite{stanley-thesis}, 
but rather is a formula developed by Stanley in \cite{stanley-thesis}.
We take it as our definition here.

\begin{definition}\label{defn:qsym-of-P}
Let $P$ be a labeled poset.
Let  $\cL(P)$ denote the \emph{Jordan-H\"older set of P}, that is, the set of all permutations in $\permsof(P)$ that are linear extensions of $P$. 
Then the \emph{quasisymmetric function of} $P$ is
\begin{equation} \label{eqn:gessels-formula}
F(P) := \sum_{\pi \in \cL(P)} L_{C(\pi)}. 
\end{equation}
\end{definition}

\begin{remark} \label{rem:relative-order}
The function $F(P)$ depends only  on the relative partial order of the labels at each covering relation of the poset and not on the absolute values of the labels themselves.
Given labeled posets $P$ defined on the set of labels $A$, and $P'$ defined on the set of labels $B$,
and a function $f:A \to B$ which is an isomorphism of their underlying unlabeled posets,
then $F(P) = F(P')$ if for every covering relation  ($y$ covers $x$) in $P$ we have
$x < y \Longleftrightarrow f(x) < f(y)$.
\end{remark}

\section{The new basis} \label{sec:new-basis}

The main goal of this section is to define our new basis (see Definition \ref{def:new-basis} below)
and to prove that it is in fact a $\bZ$-basis of $\Qsym$,
that is to say, every quasisymmetric function that can be written in terms of either the standard monomial or fundamental basis using only integer coefficients can also be written in terms of the new basis using only integer coefficients.
In Section \ref{sec:structure-constants}, we prove the positivity of the structure constants for this new basis and the grading of $\Qsym$ by composition rank.

Following the notation of \cite{billera-jia-reiner},
given \emph{unlabelled} posets $P$ and $Q$,
we denote by $P \oplus Q$ their \emph{ordinal sum}.
The set of elements of $P \oplus Q$ is the disjoint union of the elements of $P$ and $Q$.
All of the order relations of $P$ and $Q$ are retained, 
and in addition, $x <_{P \oplus Q} y$ for all $x \in P$ and $y \in Q$.

As with all the well-known bases for $\Qsym$, 
the elements of the new basis are indexed by compositions.
We denote the basis by $\{ N_\alpha \}$, where $\alpha$ ranges over all compositions.

\begin{definition} \label{def:new-basis} 
For a given composition $\alpha \neq \fatzero$, 
let $P_\alpha = A_1 \oplus \cdots \oplus A_m$ be the graded poset on $|\alpha|$ elements,
where $m = \ell(\alpha)$
and $A_i$ is an antichain on $\alpha_i$ elements.
Make $P_\alpha$ into a labeled poset by numbering the ranks in alternating fashion:
first number the odd-ranked elements $A_2, A_4, \ldots$,
followed by the even-ranked elements $A_1, A_3, \ldots$.
We define $N_\fatzero := 1$, and 
for each $\alpha \neq \fatzero$, we define $N_\alpha := F(P_\alpha)$
(see Equation  \eqref{eqn:gessels-formula}).
\end{definition}


\begin{example}
Let $\alpha = (1,2,2)$. Then
\begin{eqnarray*}
 P_{122} & = & \{3\} \oplus \{1,2\} \oplus \{4,5\}, \\
\cL(P_{122}) & = & \{ (31245),(31254), (32145),(32154) \}, \mbox{ and} \\
N_{122} & = & L_{14} + L_{131}+L_{113}+L_{1121}.
\end{eqnarray*}
\end{example}

\bigskip
\begin{definition} \label{def:rank-composition}
Given a composition $\alpha = (\alpha_1,\ldots,\alpha_k)$,
the \emph{rank} of $\alpha$, denoted by $r(\alpha)$,
is the sum of the odd-indexed parts of the composition.
That is, 
\begin{equation} \label{eqn:def-rank}
r(\alpha) := \sum_{\text{odd} \; i} \alpha_i = \alpha_1 + \alpha_3 + \alpha_5 + \cdots.
\end{equation}
We define $\cN^0_0 := \{N_\fatzero\} = \{1\}$, and for $1 \leq r \leq n$,
\begin{equation} \label{eqn:def-N-n-r}
\cN^n_r := \{ N_\alpha \suchthat |\alpha| = n \mbox{  and  } r(\alpha) = r \}.
\end{equation}
We also define the subspace $V^n_r := \spam \; \cN^n_r \subset \Qsym_n$.
If we are working over a field of coefficients for $\Qsym$, then $V^n_r$ may be viewed as a vector space, whereas if we are working over integer coefficients then we refer to the $\bZ$-span of $\cN^n_r$ and $V^n_r$ is a $\bZ$-module.
\end{definition}

\label{sec:Z-basis}
\label{sec:basis-proof}

\begin{theorem} \label{thm:Z-basis}
The set of quasisymmetric functions $\{ N_\alpha \}$, as $\alpha$ ranges over all compositions, forms a $\bZ$-basis for $\Qsym$.
\end{theorem}
\begin{proof}
We show that $\{ N_\alpha \}_{|\alpha|=n}$ forms a basis for the homogeneous component $\Qsym_n$ for each nonnegative integer $n$.
This is trivial for $n = 0$.
For the general case, we prove the existence of a unitriangular transition matrix from 
 $\{ N_\alpha \}_{|\alpha|=n}$ of $\Qsym_n$ to the fundamental basis $\{ L_\alpha \}_{|\alpha|=n}$.
 
Consider the following construction.
Given a permutation $\omega \in S_n$, let $b(\omega) \in 1\{0,1\}^{n-1}$ be the $n$-digit binary word where the digits are given by
 \[  b_i=\left\{  \begin{array}{cl}
                 1 & \mbox{if $i=1$ or } \omega(i-1) < \omega(i), \\  0 & \mbox{otherwise.}
 \end{array}  \right.
\]
 Then define  $\rho(\omega)$ to be the composition which gives the lengths of successive runs of 1's and 0's in $b(\omega)$.
 For example, if $\omega = 184356729 \in S_9$ then $b(\omega) = 110011101$, and $\rho(\omega) = 22311$.
Clearly one can determine the run-length composition $C(\omega)$ from $\rho(\omega)$
and vice versa.

Given a composition $\alpha$, let $P_\alpha$ be the labeled poset in Definition \ref{def:new-basis}, and $\cL(P_\alpha)$ its set of linear extensions.
Recall that by definition
\[ N_\alpha :=  \sum_{\pi \in \cL(P_\alpha)} L_{C(\pi)}. \]
By the nature of the labeling on $P_\alpha$, $\rho(\pi) \refines \alpha$ for all $\pi \in \cL(P_\alpha)$.
Furthermore, there is a unique element $\pi \in \cL(P_\alpha)$ such that $\rho(\pi) = \alpha$,
namely the one in which all the labels of $A_i$ are in ascending order if $i$ is odd, and in descending order if $i$ is even.
Thus if we order the rows of the transition matrix (labeled by compositions $\alpha$) and columns (labeled by compositions $\rho(\pi)$) in an arbitrary way that extends the partial refinement order $\refines$, then the resulting matrix is unitriangular, and hence $\{N_\alpha\}$ is indeed a $\bZ$-basis for $\Qsym$.
\end{proof}

\section{Additional facts regarding $F(P)$}\label{sec:F-of-p}

In this section we develop several additional facts regarding the quasisymmetric function $F(P)$ for labeled posets $P$,
including an alternative way to decompose $F(P)$ for posets,
the main idea being to partition $\cL(P)$.
These facts, especially Lemmas \ref{lem:K-1-K} through \ref{lem:chain-decomp},
are key tools for the results in following sections.

\subsection{Ordered partitions}

Consider a permutation  $\pi = (\pi_1,\ldots,\pi_n) \in \permsof(X)$, where $|X| = n$,
and a composition $\tau = (\tau_1,\ldots,\tau_k)$ with $|\tau| = n$.
We can ``chop up", or segment the one-line notation of $\pi$ from left to right into $k$ \emph{segments}, 
where each respective segment $s_i$ is a subsequence \emph{of consecutive elements} of the one-line notation of length $\tau_i$.
We call the sequence of these segments $s = (s_1,\ldots,s_k)$  a \emph{segmentation of $\pi$ of type $\tau$} (or induced by $\tau$).
Letting $t_0 = 0$ and $t_j = \sum_{i=1}^j \tau_i$ be the $j$-th partial sum of the parts of $\tau$,
every permutation $\pi \in \permsof(X)$ has a unique segmentation $s_\tau(\pi)$,
whose segments, for $1 \leq j \leq k$, are given by
\[ s_j = (\pi_{t_{j-1} + 1},\pi_{t_{j-1} + 2},\ldots,\pi_{t_j}). \]

An \emph{ordered partition} $K = (K_1,\ldots,K_k)$ of a set $X \subset \bP$ is a partitioning of $X$ into non-empty, pairwise disjoint subsets called \emph{blocks}, i.e. $X = \sqcup_{i=1}^k K_i$, where the order of the blocks matters.
Let $\tau_i = |K_i|$ for all $i$, and refer to the resulting composition $\tau(K) = (\tau_1,\ldots,\tau_k)$ as the \emph{type} of $K$.

Let $\bK(X)$ denote the set of all ordered partitions of $X$.
Every composition $\tau$ of weight $n$ induces a mapping $\cK_\tau : \permsof(X) \to \bK(X)$ as follows.
For every $\pi \in \permsof(X)$ there is a unique ordered partition $\cK_\tau(\pi)$,
each of whose blocks $K_j$ is the set of elements in the corresponding segment $s_j$ of the segmentation $s_\tau(\pi)$.

We abbreviate the inverse image $\cK_\tau^{-1}(K)$ as $\cK^{-1}(K)$ since, for a given ordered partition $K$, the type $\tau$, the set of elements $X$, and thus the permutation group $\permsof(X)$, can all be determined from $K$.
Thus for an ordered partition $K$, we have
\begin{equation} \label{eqn:K-inverse}
\cK^{-1}(K) := \{\pi \in \permsof(X) \suchthat K_{\tau(K)}(\pi) = K \}.
\end{equation}
For example, $\cK^{-1}((\{2,7\},\{5\},\{1,8\})) = \{27518, 27581, 72518, 72581\}$.

The following lemma is simply an exercise in notation, so we omit its proof.

\begin{lemma} \label{lem:K-1-K}
Let $K$ be an ordered partition with $k$ blocks. 
Let $P_K$ be the labeled poset $P_K = K_1 \oplus \cdots \oplus K_k$, where each $K_i$ is regarded as an antichain.
Then 
\[ F(P_K) = \sum_{\pi \in \cK^{-1}(K)} L_{C(\pi)}. \]
\end{lemma}

We say that an ordered partition $K = (K_1,\ldots,K_k)$ is \emph{alternating}
if for every $1 \leq i < k$ and for all $x \in K_i$ and  $y\in K_{i+1}$ we have 
$x < y$ if $i$ is even and $x > y$ if $i$ is odd.

\begin{lemma} \label{lem:alternating}
Let $K$ be an alternating ordered partition of type $\tau$.
Then \[ F(P_K) = N_\tau. \]
\end{lemma}
\begin{proof}
Each rank $K_i$ of $P_K = K_1 \oplus \cdots \oplus K_k$ (where $\ell(\tau) = k$) is an antichain.
Hence $F(P_K)$ depends only on the relative ordering of the elements between adjacent ranks $K_i$ and $K_{i+1}$ (see Remark \ref{rem:relative-order}).
Since $K$ is alternating, we can relabel its elements in each rank as we do in the construction of $P_\tau$ (as in Definition  \ref{def:new-basis}) and still maintain the same relative ordering between elements in adjacent ranks.
Thus $F(P_K) = F(P_\tau) = N_\tau$.
\end{proof}

\subsection{Unordered partitions of $X \subset \bP$}
 
Let $T = \{T_1,\ldots,T_m\}$ be an unordered partition of the set $X \subset \bP$.
We say that an ordered partition $K$ is a refinement of $T$ if $K$, considered as an unordered partition, is a refinement of $T$.
For every permutation $\pi \in \permsof(X)$,
$T$ induces a unique segmentation of $\pi$ where each segment is contained in a block of $T$ and this segmentation is least (coarsest), with respect to refinement, among all such segmentations.
Corresponding to this segmentation
there is a unique ordered partition $K_T(\pi)$, which clearly is is a refinement of $T$.
We say that $T$ \emph{induces} the ordered partition $K_T(\pi)$ on $\pi$.

\begin{example}
Let $X = [9]$, $T = \{\, \{1,4\},\, \{2,6,8,9\},\, \{3,5,7\}\, \}$, $\pi = 965412378$.
Then $K_T(\pi) = (\{6,9\},\, \{5\},\, \{1,4\},\, \{2\},\, \{3,7\},\, \{8\})$.
\end{example}

Let $P$ be a labeled poset,
and $T$ an unordered partition of $P$.
Define $\bK_{P,T}$ to be the set of induced ordered partitions
$\bK_{P,T} := \{ K_T(\pi)\;|\; \pi \in \cL(P) \}$.
We say that $T$ is \emph{antichain-inducing} if 
for every ordered partition $K \in \bK_{P,T}$, every  block $K_i$ of $K$ is an antichain in $P$.

\begin{lemma} \label{lem:chain-decomp} 
Let T be an antichain-inducing unordered partition of a labeled poset $P$.
Then
\begin{equation} \label{eqn:chain-decomp}
F(P) = \sum_{K \in \bK_{P,T}} F(P_K).
\end{equation}
We call this the \emph{decomposition of $F(P)$ with respect to $T$}.
\end{lemma}
\begin{proof}
By Lemma \ref{lem:K-1-K} it suffices to show that 
\[ \cL(P) = \bigsqcup_{K \in \bK_{P,T}} \cK^{-1}(K). \]
The  ``$\subset$"-direction  is trivial. 
Indeed, $T$ induces some ordered partition on every permutation,
and by definition $\bK_{P,T}$ includes all such partitions as permutations range over $\cL(P)$.
Also, clearly $\cK^{-1}(K) \cap \cK^{-1}(J) = \emptyset$ if $K \neq J$ since $K_T$ is a well-defined map on $\cL(P)$, 
and so the union on the right is indeed a disjoint union.

For the ``$\supset$"-direction, let $K \in \bK_{P,T}$.
By definition of $\bK_{P,T}$, there exists $\pi \in \cL(P) \cap \cK^{-1}(K)$.
Let $s = s_{\tau(K)}(\pi)$.
Since $T$ is antichain-inducing, 
the unordered set of elements $K_i$ of each segment $s_i$ is an antichain.
It follows that if we form a new permutation $\widehat\pi$ by permuting the elements of $s_i$ arbitrarily within $s_i$ (and thus within $\pi$), we must also have that $\widehat\pi \in \cL(P)$.
Since this holds true for each segment of $s$, we have
 $\cK^{-1}(K) \subset \cL(P)$.
\end{proof}

\begin{uremark}
In the extreme case where $T$ consists of all singleton sets,
$K_T(\pi)$ is the list of singleton sets in the order specified by $\pi$,
and $\cK^{-1}(K_T(\pi)) = \{ \pi \}$.
We can identify $K_T(\pi)$ with $\pi$ itself,
and similarly $\bK_{P,T}$  with $\cL(P)$,
and the lemma is then equivalent to the formula \eqref{eqn:gessels-formula}.
\end{uremark}

\subsection{Structure constants for the new basis}

Following the notation of \cite{billera-jia-reiner} and \cite{stanley-thesis},
given labeled posets $P$ and $Q$ on sets $X$ and $Y$ respectively,
we denote by $P \disjointsum Q$ any \emph{disjoint sum} of the posets, constructed as follows.
We first form the poset whose set of elements is the disjoint union of the sets of elements of $P$ and $Q$, retaining all partial order relations of the two posets but adding no new relations.
In order to ensure that all labels are distinct, we then relabel the elements in any fashion subject to the restriction that the resulting labels are all distinct and preserve the relative order of labels at all covering relations (see Remark \ref{rem:relative-order}).
While the disjoint sum of the labeled posets is not uniquely defined, all disjoint sums so constructed will have the same quasisymmetric function.
It is well-known and is easy to prove (see, for example, \cite{gessel}) that

\begin{equation}\label{eqn:disjoint-sum}
 F(P \disjointsum Q) = F(P)\cdot F(Q). 
\end{equation}
We are now in a position to prove the nonnegativity of the structure constants for our new basis.

\label{sec:structure-constants} 
\begin{theorem} \label{thm:structure-constants}
The quasisymmetric function algebra $\Qsym$ is graded by the rank of the compositions indexing the basis $\{N_\alpha\}$. 
Furthermore, the structure constants for $\{ N_\alpha \}$ are nonnegative.
That is, in the expansion
\[ N_\alpha N_\beta = \sum_\nu c^\nu_{\alpha,\beta} N_\nu, \]
all the constants $c^\nu_{\alpha,\beta}$ are nonnegative integers.
\end{theorem}

\begin{proof}
We first prove the statement regarding structure constants.
Since $N_\fatzero = 1$, the claim holds trivially if $\alpha = \fatzero$ or $\beta=\fatzero$.
Thus we assume  $\alpha = (\alpha_1,\ldots,\alpha_s) \neq \fatzero$ and $\beta = (\beta_1,\ldots,\beta_t) \neq \fatzero$.
By \eqref{eqn:disjoint-sum} we have that
\begin{equation} \label{eqn:thm-sc-first}
N_\alpha N_\beta = F(P_\alpha) F(P_\beta) = F(P_\alpha \disjointsum P_\beta).
\end{equation}
We write $P_\alpha = A_1 \oplus \cdots \oplus A_s$
and $P_\beta = B_1 \oplus \cdots \oplus B_t$, 
and identify the $A_i$ and $B_j$ subsets with their canonical inclusions in $P_\alpha \disjointsum P_\beta$.

We form a new poset $Q$ by relabeling the elements of the $A_i$ and $B_j$ subsets while maintaining their ordering relations: 
first label the even-indexed $A_i$ and $B_j$ in order, with the numbers from $[m]$,
where $m = |\alpha| + |\beta| - r(\alpha) - r(\beta)$
and $r(\alpha)$ is the rank function from Definition \ref{def:rank-composition},
then label the odd-indexed $A_i$ and $B_j$ in order, with the numbers from $\{m+1,\ldots,|\alpha| + |\beta| \}$.
Since $F(P_\alpha \disjointsum P_\beta)$ depends only on the relative ordering of elements between adjacent ranks $A_i$ and $A_{i+1}$ for $1 \leq i < s$ 
and between adjacent ranks $B_j$ and $B_{j+1}$ for $1 \leq j < t$,
we have 
\begin{equation}
F(P_\alpha \disjointsum P_\beta) = F(Q).
\end{equation}
We consider the unordered partition $T = \{T_1,T_2\}$ of $Q$ given by 
 \[ T_1 = \left(\bigcup_{odd\; i} A_i \right) \bigcup \left(\bigcup_{odd\; i} B_i \right),
\quad \mbox{ and } \quad
T_2 = \left(\bigcup_{even\; i} A_i \right) \bigcup \left(\bigcup_{even\; i} B_i \right). \]
Note that $T$ is antichain-inducing, so we may apply Lemma \ref{lem:chain-decomp}:
\begin{equation}
F(Q) = \sum_{K \in \bK_{Q,T}} F(P_K).
\end{equation}
On the other hand, the labeling of $Q$ implies that every ordered partition $K \in \bK_{Q,T}$ is alternating,
so applying Lemma \ref{lem:alternating}, we have
\begin{equation}\label{eqn:thm-sc-last}
F(Q) = \sum_{K \in \bK_{Q,T}} N_{\tau(K)}.
\end{equation}
Combining Equations \eqref{eqn:thm-sc-first} -- \eqref{eqn:thm-sc-last} yields the positivity claim.
In particular, \[ c^\nu_{\alpha,\beta} = |\{ K \in \bK_{Q,T} \suchthat \tau(K) = \nu \}|. \]

To prove the statement regarding the grading of $\Qsym$ by composition rank,
we simply note that for every $K \in \bK_{Q,T}$, we have \[ r(\tau(K)) = |T_1| = r(\alpha) + r(\beta). \]
\end{proof}

\section{Matroids} \label{sec:matroids}

This section begins the second part of the paper.
Here we review some of the concepts, terminology, and results from  \cite{billera-jia-reiner},
and prove our claims regarding the quasisymmetric functions of matroids vis-a-vis our new basis.
For general background in matroid theory we refer the reader to standard texts such as Oxley's \cite{oxley}.
We review several of the terms here.

The \emph{direct sum} of matroids $M_1$ and $M_2$, denoted $M_1 \oplus M_2$, has as its ground set the disjoint union $E(M_1 \oplus M_2) = E(M_1) \sqcup E(M_2)$, and as its bases 
\[ \cB(M_1 \oplus M_2) = \{ B_1 \sqcup B_2 \suchthat B_1 \in \cB(M_1), B_2 \in \cB(M_2) \}. \]
A \emph{circuit} is a minimal dependent set.
If we declare two elements of a matroid to be equivalent if and only if they are both contained in some circuit, then the equivalence classes of elements are the \emph{components} of the matroid.
We say that the matroid is \emph{connected} if it has only one component, and \emph{disconnected} otherwise.
A matroid is the direct sum of its components.

\subsection{The quasisymmetric function of a matroid}

\BJR \cite{billera-jia-reiner} describe an invariant for isomorphism classes of matroids in the form of a quasisymmetric function.  
Rather than give the definition from   \cite{billera-jia-reiner}, we describe it in terms of a formula which is shown in  \cite{billera-jia-reiner}  to be equivalent to the definition.

Fix a matroid $M$, one of its bases $B \in \cB(M)$, and let $B^c = E(M) - B$ (the \emph{cobase} of $B$).
Define the poset $P_B$ on the ground set $E(M)$ 
where $e <_{P_B} e'$ if and only if $e \in B$, $e' \in B^c $, and $(B-e)\cup\{e'\} \in \cB(M)$.
That is, $e <_{P_B} e'$ if and only if swapping $e'$ for $e$ in $B$ yields another base in $M$.
Thus the Hasse diagram of $P_B$ is a bipartite graph 
in which the elements of $B$ are minimal elements of the poset and the elements of $B^c$ are maximal elements.
Note that if $M$ has no loops, then in the Hasse diagram of $P_B$, 
every element in $B^c$ has positive vertex degree.
We say that a labeled poset is \emph{strictly labeled} if for all $x,y \in P$ we have that $x <_P y$ implies $x > y$.
Similarly,  a labeled poset is  \emph{naturally labeled} if for all $x,y \in P$, $x <_P y$ implies $x < y$.
We apply a strict labeling to $P_B$ (any will do).
The quasisymmetric function $F(M)$ associated with $M$ can be written as
\begin{equation} \label{eqn:base-poset}
F(M) = \sum_{B \in \cB(M)} F(P_B),
\end{equation}
where $F(P_B)$ is the quasisymmetric function of the strictly labeled poset $P_B$ as defined in Definition \ref{defn:qsym-of-P}. 

It was shown in \cite{billera-jia-reiner} that the mapping  $F: \Mat \to \Qsym$ is in fact a morphism of combinatorial Hopf algebras, with a suitable choice of character on the algebra $\Mat$.
Here $\Mat$ is the Hopf algebra of matroids introduced by Schmitt \cite{schmitt-I} and studied by Crapo and Schmitt 
\cite{crapo-schmitt-4}, \cite{crapo-schmitt-1}, \cite{crapo-schmitt-2}, \cite{crapo-schmitt-3}.
The matroid algebra $\Mat$ has as its basis elements isomorphism classes of matroids.
The product of two basis elements $[M_1]$ and $[M_2]$ in the algebra is given by
$[M_1] \cdot [M_2] := [M_1 \oplus M_2]$, where $M_1 \oplus M_2$ denotes the \emph{direct sum} of matroids.
Comultiplication in $\Mat$ is given by 
$\Delta([M]) := \sum_{A \subset E(M)} [M|_A] \otimes [M\setminus A]$,
where $M|_A$ is the restriction of $M$ to $A$,
and $M\setminus A$ is the contraction of $M$ by $A$.
Under the morphism $F$ we have that
\[ F(M_1 \oplus M_2) = F(M_1) \cdot F(M_2). \]
\BJR  also show that the mapping $F$, while not surjective over integer coefficients,  is surjective over rational coefficients.

\bigskip
\label{sec:loopless}
The mapping  $F: \Mat \to \Qsym$  does not distinguish between loops and coloops. 
Indeed,  let $M'$ extend the matroid $M$ by adding a loop $\ell$,
i.e. $M' = M \oplus \{\ell\}$,
and  let $M''$ extend the matroid $M$ by adding a coloop $c$,
i.e. $M'' = M \oplus \{c\}$.
Then
\begin{equation} \label{eqn:loop-coloop} 
F(M') = F(M'') = F(M) \cdot L_1.
\end{equation}
Here $L_1$ is the fundamental basis function indexed by the composition $(1)$,
which is the  elementary symmetric function $e_1(\xx)$.
Define an equivalence relation $\sim$ on isomorphism classes of matroids by 
$[M_1] \sim [M_2]$ if and only if one can obtain $M_1$ from $M_2$ by changing some number of loops to coloops or vice versa.
Then by Equation \eqref{eqn:loop-coloop} the mapping  $\Mat \to \Qsym$  
factors through the quotient
\[ \Mat \rightarrow \Mat/\!\!\sim\; \rightarrow \Qsym. \]
Accordingly, throughout most of our paper, we assume that, unless otherwise specified, 
our matroids have no loops;
that is, out of each equivalence class in $\Mat/\!\!\sim$ we select the representative that has no loops when considering their images in $\Qsym$.

\subsection{Expanding $F(M)$ in the $\{N_\alpha\}$ basis}

Recall from Definition \ref{def:rank-composition} that 
$\cN^n_r = \{ N_\alpha \suchthat |\alpha| = n, r(\alpha) = r \}$ 
and $V^n_r = \spam\: \cN^n_r$.
In this subsection we may take our coefficient ring to be $\bZ$ if we wish.

\begin{lemma} \label{lemma:rank1}
Let $P$ be a strictly labeled poset on $n$ elements of rank at most one and with $r$ minimal elements.
Then $F(P) \in V^n_r$.
Moreover the expansion of $F(P)$ in terms of the basis elements $\cN^n_r$ has only nonnegative integer coefficients.
\end{lemma}
\begin{proof}
If $P$ has rank 0, then $P$ is an antichain and so $r = n$.
Thus by labelling $P$ with the elements of $[n]$ and taking $\alpha = (n)$, we have
\[ F(P) = F(P_\alpha) = N_\alpha = N_{(n)} \in V^n_n. \]

Otherwise $P$ has rank 1 and is not an antichain.
Let   $T = \{T_1,T_2\}$ be the unordered partition of $P$ 
in which $T_1$ comprises the $r$ minimal elements of $P$, 
and $T_2$  the remaining elements.  
Note that some elements may be both minimal and maximal, and these will be placed in $T_1$.
Thus every element in $T_2$ in the Hasse diagram of $P$ has positive vertex degree.
Since we are interested in computing $F(P)$, we may assume without loss of generality
that $P$ has a strict labeling which labels the elements of 
$T_1 = \{ n, n-1, \ldots, n-r+1 \}$ in arbitrary fashion,
and the elements of $T_2 = \{ 1, 2, \ldots, n-r \}$ in arbitrary fashion.
Since $T_1$ and $T_2$ are themselves antichains, 
$T$ is antichain-inducing, so by Lemma  \ref{lem:chain-decomp}:
\[ F(P) = \sum_{K \in \bK_{P,T}} F(P_K). \]
Moreover, by the choice of labeling and the fact that every element in $T_2$ has positive vertex degree,
we have that every $K \in \bK_{P,T}$ is alternating.
Lemma \ref{lem:alternating} then implies
\[ F(P) = \sum_{K \in \bK_{P,T}} N_{\tau(K)}. \]
We also have that $|\tau(K)| = n$ and $r(\tau(K)) = \sum_{odd \; i} \tau_i = |T_1| = r$,
thus $N_{\tau(K)} \in V^n_r$ for every $K \in  \bK_{P,T}$.
Hence $F(P) \in V^n_r$ as claimed.
\end{proof}

\label{sec:matroid-friendly} 
\begin{theorem} \label{thm:FM-membership}
Let $M$ be a loopless matroid of rank $r$ on $n$ elements.
Then $F(M) \in V^n_r$.  
Moreover the expansion of $F(M)$ in terms of the basis elements $\cN^n_r$ has only nonnegative integer coefficients.
\end{theorem}
\begin{proof}
For a loopless matroid $M$, for every base $B \in \cB(M)$,  
the base poset $P_B$ has $r$ minimal elements out of a total of $n$ elements, 
and  rank at most one. 
The assertion then follows by Lemma \ref{lemma:rank1}, 
and the formula in Equation \eqref{eqn:base-poset}.
\end{proof}

\label{sec:N-coeffs}
We make a few observations here about the coefficients in the expansion of the quasisymmetric function $F(M)$ of a matroid $M$ in terms of our new basis.
Given a quasisymmetric function $q$, define 
\[ \supp(q) = \left\{ \alpha \suchthat m_\alpha \neq 0 \mbox{ in the expansion } q = \sum_\alpha m_\alpha N_\alpha \right\}. \]
We first note that for a (loopless) matroid $M$ of rank $r$ on $n$ elements,
the coefficient $m_{(r,n-r)}$ of $N_\alpha$, where $\alpha = (r,n-r)$,
is equal to the number of bases of $M$.

\begin{example} \label{exm:uniform}
Let $M = U_{r,n}$ be the uniform matroid of rank $r$ on $n$ elements.
By definition, its bases are all the $r$-subsets of the ground set $E(M)$, i.e. $\cB(M) ={E(M) \choose r}$.
Then every $P_B$ has a complete bipartite $K_{r,n-r}$ graph for its Hasse diagram.
Therefore $F(U_{r,n}) = {n \choose r} N_{r,n-r}$, and $\supp(F(M)) = \{ (r,n-r) \}$.
\end{example}

Thus the values of $m_{(r,n-r)}$ (how small) and $|\supp(F(M))|$ (how large) 
are, to some extent, measures of the degree to which $M$ fails to be uniform.

\bigskip
The coefficient $m_\alpha$ where $\alpha = (r-1,1,1,n-r-1)$ also has a combinatorial interpretation.
There is such an $N_\alpha$ term for every edge ``missing" from the Hasse diagram of a base poset $P_B$ as compared to the complete bipartite graph $K_{r,n-r}$.
In terms of matroid base polytopes, which are discussed in \ref{sec:polytopes},
the polytope $Q(U_{r,n-r})$ contains all possible vertices, namely ${E(M) \choose r}$,
while the base polytope for a different matroid $M$ of same rank and ground set size 
has only a subset of them, namely $\cB(M)$.
The coefficient $m_\alpha$ is the number of edges in the 1-skeleton of $Q(U_{r,n-r})$
between the set of vertices $\cB(M)$ and its complement ${E(M) \choose r} - \cB(M)$.

\begin{lemma} \label{lem:num-coloops}
Let $M$ be a matroid, possibly containing loops.
Then the total number $c$ of loops and coloops is given by
\begin{equation} \label{eqn:num-coloops}
c = \max_{\alpha \in \supp(F(M)), \atop \ell(\alpha) \; \text{odd}}  \dot{\alpha},
\end{equation}
where $\dot{\alpha}$ denotes the last part of the composition $\alpha$.
\end{lemma}
\begin{proof}
Since the morphism $F:\Mat \to \Qsym$ factors through loop-coloop equivalence,
$F(M) = F(M')$ where $M'$ is obtained from $M$ by replacing all loops of $M$ with coloops.
Since $M$ and $M'$ both have the same total number of loops and coloops,
without loss of generality, we assume that $M$ has no loops.

Consider a typical strictly labeled base poset $P_B$ of M, 
and antichain inducing partition $\{B,B^c\}$ as in the proof of Lemma \ref{lemma:rank1}.
We have $\alpha \in \supp(F(M))$ if and only if
there is an induced ordered partition $K$ of $P_B$ of type $\alpha$.
Now $\ell(\alpha)$ is odd if and only if the last block of $K$ is a subset of $B$, 
and in this case the elements in this block must be coloops. Thus $c \geq \dot{\alpha}$.
Conversely, there always exists an induced ordered partition of the poset, say of type $\alpha$, which has all of the coloops of $M$ in the last block, i.e. $c = \dot{\alpha}$, and this ordered partition will have odd length.  The result follows.
\end{proof}

%
%
\section{Matroid base polytopes} \label{sec:rank2}

In this section we recall  the context of  \cite[Section 7] {billera-jia-reiner} regarding
the relationship between decompositions of the quasisymmetric function associated to a matroid
and decompositions of its matroid base polytope.
In Subsection \ref{sec:geom} we develop a formula for the quasisymmetric function of a loopless rank two matroid in terms of the new basis, and apply it to address \cite[Question 7.12] {billera-jia-reiner} and \cite[Question 7.10] {billera-jia-reiner}.

\subsection{Matroid base polytopes and their decompositions} \label{sec:polytopes}

The motivating context is to study the decompositions of the \emph{matroid base polytope} $Q(M)$ of a matroid $M$.
This topic arises in the work of Lafforgue \cite{lafforgue-1}, \cite{lafforgue-2}, Kapranov \cite[\S 1.2 -- 1.4]{kapranov}, and can be found in the work of Speyer \cite{speyer}.

If $M$ is a matroid with $|E(M)|=n$, we define the  \emph{matroid base polytope} $Q(M)$ 
by identifying $E(M)$ with the set of standard basis vectors $\{e_i\}_{i=1}^n$ of $\bR^n$ and declaring
\[ Q(M) := \conv \left\{  \sum_{e_i \in B} e_i \suchthat B \in \cB(M) \right\}, \]
where $\cB(M)$ is the set of bases of $M$.
Useful facts about matroid base polytopes (see \cite{gelfand-serganova}), which we quote without proof, are:
\begin{enumerate}
\item 
If $M$ has rank $r$, then $Q(M)$ lies in the hyperplane $\{ \xx \in \bR^n \suchthat \sum_i x_i = r \}$.
\item
There is an edge in $Q(M)$ between vertices (bases) $B_1$ and $B_2$ if and only if there exist a pair of elements $e_i \in B_1$ and $e_j \in B_2$ such that $B_2 = (B_1 - \{e_i\}) \cup \{e_j\}$.
\item
Each face of a matroid base polytope is in turn the base polytope of some matroid.
\item   \label{Qfact:dimension}
The dimension of $Q(M)$ is $|E(M)| - s(M)$, where $s(M)$ is the number of connected components of $M$.
\end{enumerate}
\BJR define a \emph{matroid base polytope decomposition} of $Q(M)$ to be a decomposition 
\begin{equation}  \label{eqn:decomp-QMi}
Q(M) = \bigcup_{i=1}^t Q(M_i), 
\end{equation}
where each $Q(M_i)$ is also a matroid base polytope for some matroid, and for each $i \neq j$, the intersection $Q(M_i) \cap Q(M_j) = Q(M_i \cap M_j)$ is a face of both $Q(M_i)$ and $Q(M_j)$.
They call such 
a decomposition a \emph{split} if $t = 2$.

\BJR show that the mapping $F :\Mat \to \Qsym$ behaves as a \emph{valuation}  on matroid base polytopes. (See \cite{barvinok} for a discussion of valuations.)
This implies that, given a matroid base polytope decomposition as in Equation \eqref{eqn:decomp-QMi},
$F(M)$ can be expressed in terms of the set of $F(M_j)$ in an inclusion-exclusion fashion, 
where the $M_j$ are the matroids of the faces of the constituent polytopes in the decomposition.
For example, given a split $Q(M) = Q(M_1) \cup Q(M_2)$, we have
$F(M) = F(M_1) + F(M_2) - F(M_1 \cap M_2)$,
where $Q(M_1 \cap M_2)$ is necessarily a lower-dimensional face.
Hence by Fact \ref{Qfact:dimension} above, the matroid $M_1 \cap M_2$ is disconnected,
and so $F(M_1 \cap M_2)$ can be expressed as a product.
If we let $\frakm := \bigoplus_{d \geq 1} \Qsym_d$ be the maximal ideal in the ring $\Qsym$,
then $F(M_1 \cap M_2) \in \frakm^2$.
Therefore in the quotient space $\Qsym/\frakm^2$, we have $\overline{F(M)} = \overline{F(M_1)} + \overline{F(M_2)} $.
In general, given a matroid base polytope decomposition as in Equation \eqref{eqn:decomp-QMi},
 there is an algebraic decomposition modulo $\frakm^2$
\begin{equation} \label{eqn:decomp-oFMi}
\overline{F(M)} = \sum_i \overline{F(M_i)}.
\end{equation}
One of the open questions raised by \BJR \cite{billera-jia-reiner} is under what conditions the converse may hold;
given a collection of matroids satisfying \eqref{eqn:decomp-oFMi}, what additional conditions are sufficient to conclude \eqref{eqn:decomp-QMi}?
\begin{note}
So far we have ignored the distinction between the isomorphism class of a matroid on the one hand and a specific instance of that class on a given ground set on the other, since the quasisymmetric function of a matroid is invariant on the elements of the same isomorphism class.
When discussing the existence of matroid base polytope decompositions, it is sometimes necessary to draw a distinction between the notions, as is done in the statement of Theorem \ref{thm:geom} below.
When this precision is necessary, we use the usual bracket notation $[M]$ to denote the isomorphism class of the matroid $M$.
\end{note}

Given a ground set size $n$, the converse question is trivial for rank $0$ and $1$, and by matroid duality, for rank $n$ and $n-1$.
One necessary condition they point out is that a specific set of matroids on a common ground set satisfying \eqref{eqn:decomp-QMi} must at least satisfy the condition $\cB(M_i) \subset \cB(M)$ for all $i$, in which case they say that \eqref{eqn:decomp-oFMi} is a \emph{weak image decomposition} 
and that $\overline{F(M)}$ is \emph{weak image decomposable}.
They specifically ask,
\begin{quote}
\textbf{\cite[Question 7.12]{billera-jia-reiner}}
Does $\overline{F(M)}$ being weak image decomposable in $\Qsym/\frakm^2$ imply that $Q(M)$ is decomposable?
\end{quote}
So far, general sufficient conditions are not known beyond the trivial ranks listed above.
We claim that the converse (\eqref{eqn:decomp-oFMi} $\Rightarrow$ \eqref{eqn:decomp-QMi}) holds quite generally for rank two matroids, as shown in Section \ref{sec:geom}.
By matroid duality, the converse also holds for matroids of corank two.

\bigskip
As discussed in \cite{billera-jia-reiner}, the loopless rank two matroids are indexed, up to isomorphism, by partitions having two or more parts, where there are as many parts as there are parellelism classes of elements in the matroid and the parts of the partition give the respective cardinalities of these classes.
For this section we write $M_\lambda$ to denote the loopless rank two matroid indexed by the partition $\lambda$.
More generally, given a composition $\alpha$, define $M_\alpha = M_\lambda$ where $\lambda$ is the decreasing rearrangement of the parts of $\alpha$.

Kapranov \cite[\S 1.3]{kapranov} gives a description of all decompositions of rank two matroid base polytopes.
He shows  \cite[Lemma 1.3.14]{kapranov}  that in rank two, all matroid base polytope decompositions arise from hyperplane splits.
We provide some description here of the geometric situation, in our own words.
Given the composition $\lambda = (\lambda_1,\ldots,\lambda_m)$, with $|\lambda| = n$,
set $t_0 = 0$ and for $0 \leq k \leq m$ set $t_k = \sum_{i=1}^k \lambda_i$ be the $k$-th partial sum.
Then the vertices of $Q(M_\lambda)$ are precisely those $0/1$-lattice points $v$ lying in the hyperplane $H = \{ \xx \in \bR^n \suchthat \sum_i x_i = 2 \}$
subject to the restriction that \[ \sum_{i=t_k+1}^{t_{k+1}} v_i \leq 1 \] for all $0 \leq k < m$.
If $\ell(\lambda) = 2$, then 
\[ M_\lambda = M_{(\lambda_1,\lambda_2)} = U_{1,\lambda_1} \oplus U_{1,\lambda_2}, \]
where $U_{1,n}$ is the uniform matroid of rank 1 on $n$ elements.
It follows that if $\ell(\lambda) = 2$, then $\dim Q(M_\lambda) = n-2$ and $F(M_\lambda) \in \frakm^2$.

Supposing that $\ell(\lambda) =m  > 3$, choose index $j$ such that $1 < j < m-1$.
Let $a = t_j$ and $b = n - t_j$, and define compositions
$\mu = (a,b)$, $\alpha = (a,\lambda_{j+1},\ldots,\lambda_m)$, and $\beta = (\lambda_1,\ldots,\lambda_j,b)$, all of which have weight $n$.
Consider the hyperplane $H' = \{ \xx \in \bR^n \suchthat \sum_{i=1}^{t_j} x_i = 1 \}$.
Then $H' \cap Q(M_\lambda) = Q(M_\mu)$, 
giving us a hyperplane split $Q(M_\lambda) = Q(M_\alpha) \cup Q(M_\beta)$.
It follows from the above that $F(M_\lambda) = F(M_\alpha) + F(M_\beta) - F(M_\mu)$,
and  $\overline{F(M_\lambda)} = \overline{F(M_\alpha)} + \overline{F(M_\beta)} $.
We can summarize this in the following proposition.
The relations given in the proposition remain true even if $\lambda$ has only two or three parts,
but in that case the resulting relations are  trivial.

\begin{proposition} \label{prop:rank2-relations}
Let $\lambda = (\lambda_1,\ldots,\lambda_t)$ be a composition with at least two parts.
Let $1 \leq s < t$, $a = \sum_{i=1}^s \lambda_i$, and $b = \sum_{i=s+1}^t \lambda_i$.
Consider compositions 
$\alpha = (a,\lambda_{s+1},\ldots,\lambda_t)$, 
$\beta = (\lambda_{1},\ldots,\lambda_s,b)$,
and $\mu = (a,b)$.
We then have \[ F(M_\lambda) = F(M_\alpha) + F(M_\beta) - F(M_\mu), \]
and modulo $\frakm^2$, \[ \overline{F(M_\lambda)}  = \overline{F(M_\alpha)} + \overline{F(M_\beta)}. \]
Moreover there is a split of matroid base polytopes
\[ Q(M_\lambda) = Q(M_\alpha) \cup Q(M_\beta). \]
\end{proposition}
\noindent
The splitting process can be repeated on the constituent matroid base polytopes until we have decomposed $Q(M_\lambda)$ into the union of matroid base polytopes of type $Q(M_\alpha)$ where $\ell(\alpha) = 3$.
Consequently, modulo $\frakm^2$, $\overline{F(M_\lambda)}$ can be written as a positive sum
\[ \overline{F(M_\lambda)} = \sum_i \overline{F(M_i)}, \]
where each $M_i$ is a loopless rank 2 matroid indexed by a partition of length 3.

In this setting, \BJR,  pose the following question:
\begin{quote}
\textbf{\cite[Question 7.10]{billera-jia-reiner}}
Fix $n$ and consider the semigroup generated by $\overline{F(M)}$ within $\Qsym_n/\frakm^2$ as one ranges over all matroids $M$ of rank 2 on $n$ elements.
Is the Hilbert basis for this semigroup indexed by those $M$ for which $\lambda(M)$ has exactly 3 parts?
\end{quote}

By repeated application of Proposition \ref{prop:rank2-relations}, the set
$\{\overline{F(M_\lambda)} \suchthat \ell(\lambda) = 3 \}$
generates the semigroup in question,
so the point of the question is whether this generating set is minimal,
and whether distinct indices yield distinct functions.
We prove that this is the case as a corollary of Theorem  \ref{thm:geom}.

\subsection{Results for rank two matroids}  \label{sec:geom}

In this section, we prove that the morphism $F:\Mat \to \Qsym$ distinguishes isomorphism classes of rank two matroids and that decomposability of $F(M)$ for a rank two matroid $M$ implies decomposability of $Q(M)$, as stated in the following theorem.

\begin{theorem} \label{thm:geom}
Let $\lambda \vdash n$ with $\ell(\lambda) \geq 3$, and let $J$ be a multiset of partitions of $n$, all of length three or more, such that
\begin{equation}  
\overline{F([M_\lambda])} = \sum_{\mu \in J} \overline{F([M_\mu])},
\end{equation}
where $[M_\tau]$ denotes the isomorphism class of (loopless) rank two matroids on $n$ elements indexed by the partition $\tau$.
Then, taking the set of standard basis vectors of $\bR^n$ as the common ground set,
there exists a collection of representative matroids on this ground set,
$M_\lambda \in [M_\lambda]$ and $M_\mu \in [M_\mu]$ for all $\mu \in J$
which form a decomposition of matroid base polytopes
\begin{equation}  
Q(M_\lambda) = \bigcup_{\mu \in J} Q(M_\mu).
\end{equation}
\end{theorem}

Before the main proof of this theorem, we establish some preliminary results.
We begin by developing a formula for $F(M_\lambda)$ in terms of the new basis $\{N_\alpha\}$. 
We define the following quasisymmetric functions in 
\[ V^n_2 = \spam \{ N_\alpha \suchthat |\alpha|=n, r(\alpha)=2 \}. \]
For all $1 \leq k \leq n-1$ let
\begin{equation} \label{eqn:T-vecs}
T^n_k := \frac{1}{2}N_{(2,n-2)} + \sum_{j \ge 1}  {k - 1 \choose j}  N_{(1,j,1,n-2-j)},
\end{equation}
where we understand $N_{(1,j,1,n-2-j)}$ to be $N_{(1,n-2,1)}$ when $j = n-2$.
We also define quasisymmetric functions 
\[ U^n_k := k(n-k) T^n_k. \]
Note that each of the sets $\{T^n_k \}$ and $\{U^n_k \}$ forms a basis for the subspace $V^n_2$,
where we consider $\Qsym$ to have rational coefficients.

\begin{lemma} \label{lem:rank2-formula}
Let $M_\lambda$ be the rank two matroid on $n$ elements indexed by the partition $\lambda = (\lambda_1,\ldots,\lambda_m)$. 
Then
\begin{equation} \label{eqn:rank2-formula}
F(M_\lambda) =  \sum_{i=1}^m U^n_{\lambda_i}.
\end{equation}
\end{lemma}
\begin{proof}
We write $c(\lambda_i)$ to denote the parallelism class of elements in $M_\lambda$ corresponding to the part $\lambda_i$. 
A typical base $B \in \cB(M_\lambda)$ is $B = \{e_i,e_j\}$, where $e_i \in c(\lambda_i)$ and $e_j \in c(\lambda_j)$ are in distinct parallelism classes.
The Hasse diagram of $P_B$ has two minimal elements, $e_i$ and $e_j$.
There are edges from $e_i$ to all elements of the cobase $B^c = E(M_\lambda) - B$ except for the $\lambda_j -1$ elements which are in the same parallelism class $c(\lambda_j)$ as $e_j$.
Similarly,
there are edges from $e_j$ to all elements of the cobase except for the $\lambda_i -1$ elements which are in the same parallelism class $c(\lambda_i)$ as $e_i$.

We can analyze $F(P_B)$ as in the proof of Lemma \ref{lemma:rank1} by applying a strict labeling
$\gamma : E(M_\lambda) \to [n]$
such that $\gamma(e_i) = n$, $\gamma(e_j) = n-1$, and the cobase elements are arbitrarily labeled with $\{1,2,\ldots,n-2\}$.
We take $T = \{B,B^c\}$ to be our antichain-inducing partition of $(P_B,\gamma)$.
There is one induced ordered partition (of $[n]$)  of type $(2,n-2)$,
namely $K = (B, B^c)$, classifying one set of permutations in $\cL(P_B,\gamma)$, 
and thus contributing one $N_{(2,n-2)}$ term to the expansion of $F(P_B)$.
For each $1 \le k < \lambda_j$, and for each $k$-set  $A \subset B^c \cap c(\lambda_j)$,
there is an induced ordered partition $K = ( \{e_j\},A,\{e_i\},B^c-A )$ of type $(1,k,1,n-2-k)$
contributing a term $N_{(1,k,1,n-2-k)}$ to the expansion.
Thus there are ${ \lambda_j -1 \choose k }$ such terms $N_{(1,k,1,n-2-k)}$ corresponding to 
ordered partitions $K$ of type $(1,k,1,n-2-k)$ with $K_1 = \{e_j\}$.
Likewise there are ${ \lambda_i -1 \choose k }$ such terms $N_{(1,k,1,n-2-k)}$ corresponding to 
ordered partitions $K$ of type $(1,k,1,n-2-k)$ with $K_1 = \{e_i\}$.
All the $N_\alpha \in \cN^n_2$ are of one of these types, 
and we know that the terms of $F(P_B)$ must lie in $V^n_2$,
so these are the only types appearing in the expansion for $F(P_B)$.
There can be no other terms than these due to the order relations in $P_B$.
Thus 
\begin{equation}  
F(P_B) =  N_{(2,n-2)} + 
	\sum_{k \ge 1} \left(  {\lambda_i - 1 \choose k} + {\lambda_j - 1 \choose k} \right) N_{(1,k,1,n-2-k)} .
\end{equation}
Using Equation \eqref{eqn:T-vecs}, we can rewrite this as
\[ F(P_B) =  T^n_{\lambda_i} + T^n_{\lambda_j}. \]
Finally, there are $\lambda_i \lambda_j$ such bases $B \in c(\lambda_i) \times c(\lambda_j)$.
Summing over all pairs of parallelism classes of the matroid yields the formula
\[ F(M_\lambda) =   \sum_{i=1}^m \lambda_i(n-\lambda_i) T^n_{\lambda_i}
      = \sum_{i=1}^m  U^n_{\lambda_i}. 
\]
\end{proof}

Next we develop a similar formula for $\overline{F(M_\lambda)}$ in $\Qsym_n/\frakm^2$.
Our starting point is the following corollary.

\begin{corollary} \label{cor:rank2products}
Let $a$ and $b$ be positive integers such that $a + b = n$.
Then \[ ab \cdot N_{(1,a-1)} \cdot N_{(1,b-1)} = U^n_a + U^n_b. \]
\end{corollary}
\begin{proof}
Let $\lambda = (a,b)$. 
Then $M_\lambda = U_{1,a} \oplus U_{1,b}$, 
where $U_{1,m}$ is the uniform matroid of rank one on $m$ elements.
As discussed in Example  \ref{exm:uniform}, $F(U_{1,m}) = m N_{(1,m-1)}$.
Therefore by the Hopf algebra morphism, we have
\[ F(M_\lambda) = F(U_{1,a})\cdot F(U_{1,b}) = a N_{(1,a-1)} \cdot b N_{(1,b-1)}. \]
On the other hand, by Lemma \ref{lem:rank2-formula} we have $F(M_\lambda) =   U^n_a + U^n_b$.
Equating right hand sides  yields the desired formula.
\end{proof}

Since $\Qsym$ with respect to its product structure is graded by composition rank as well as degree,
the vector subspace $V^n_2 \cap \frakm^2$ is spanned by the vectors 
\[ \{ N_{(1,a-1)} \cdot N_{(1,b-1)} \suchthat  a + b = n \}. \]
Thus a basis for $V^n_2 \cap \frakm^2$ is $\{ U^n_k + U^n_{n-k} \suchthat 1 \leq k \leq \frac{n}{2} \}$.
For expressing our formula for $\overline{F(M_\lambda)}$, 
we find it convenient to define vectors $\overline{U^n_k}$ as follows:
\begin{equation} \label{eqn:def-Ubar}
\overline{U^n_k} =\left\{ \begin{array}{l@{\quad}l}
	U^n_k  &  \mbox{ if } k < \frac{n}{2}, \\
	0 &  \mbox{ if } k = \frac{n}{2}, \\
	-U^n_{n-k}  &  \mbox{ if } k > \frac{n}{2} .
\end{array} \right. 
\end{equation}
Thus the set $\{\overline{U^n_k} \suchthat 1 \leq k < \frac{n}{2} \}$ forms a basis (over rational coefficients) for $V^n_2/\frakm^2$.
We have the immediate corollary of Lemma  \ref{lem:rank2-formula}:

\begin{corollary}  \label{cor:rank2-formula-modm2}
Let $M_\lambda$ be the rank two matroid on $n$ elements indexed by the partition $\lambda = (\lambda_1,\ldots,\lambda_m)$. 
Then
\begin{equation} \label{eqn:rank2-formula-modm2}
\overline{F(M_\lambda)} =  \sum_{i=1}^m \overline{U^n_{\lambda_i}}.
\end{equation}
\end{corollary}

\noindent
The next proposition provides a necessary step for the main result, but may be of interest in its own right.

\begin{proposition} \label{prop:F-inject}
Let $\cM^2$ be the set of matroid isomorphism classes (including those with loops) of rank two matroids.
Let $\Mat_c$ be the vector subspace of $\Mat$ spanned by the isomorphism classes of connected matroids,
and let $\cM^2_c$ be the set of matroid isomorphism classes of connected rank two matroids.
Then the algebra morphism $F:\Mat \to\Qsym$ is injective when restricted to $\cM^2$, 
and the induced quotient map of vector spaces $\overline{F}:\Mat_c \to \Qsym/\frakm^2$ is injective when restricted to $\cM^2_c$.
\end{proposition}
\begin{proof}
We show that we can recover the isomorphism class of the matroid from its respective function.
Suppose we are given $F(M)$ for a rank two matroid $M$.
We know that $F(M)$ is a non-zero homogeneous function of degree $n = |E(M)|$, and so we recover the size of the ground set.  Clearly, $n \geq 2$.

It is possible that $M$ may have loops or coloops.
By Lemma \ref{lem:num-coloops} we can recover the total number $s$ of loops and coloops of $M$ from $F(M)$ by Equation  \eqref{eqn:num-coloops}.  
If $s = n$, then $M$ consists of two coloops and $n-2$ loops.
Otherwise $s \leq n-2$, and we may factor $F(M)$ as 
\begin{equation*}
F(M) = N_{(s)} \cdot F(M'),
\end{equation*}
where $(s)$ is the one-part composition of $s$,
and $M'$ is the matroid obtained from $M$ by removing all loops and coloops.
If now $F(M') \in V^{n-s}_1$, we have $M' \cong U_{1,n-s}$ and $M$ has one coloop and $s-1$ loops.
Otherwise $M$ has $s$ loops, no coloops, $F(M') \in V^{n-s}_2$, 
and $M'$ is a loopless rank two matroid on $n-s$ elements.

So now without loss of generality, we assume that $M$ has no loops or coloops
and thus is isomorphic to $M_\lambda$ for some $\lambda \vdash n$.
We expand $F(M)$ as
\begin{equation}
F(M) = \sum_{k=1}^{n-1} t_k U^n_k.
\end{equation}
This expansion can be determined since the set of $\{U^n_k\}$ form a basis of $V^n_2$.
Per Lemma \ref{lem:rank2-formula}, for each $k$, the coefficient $t_k$ is  the number of parts of $\lambda$ that are equal to $k$,
and so we recover $\lambda$ from $F(M_\lambda)$.

The argument for recovering $M$ from $\overline{F(M)}$ for a connected rank two matroid $M$
is similar.
Since $M$ is connected, it has no loops or coloops, and so again $M$ is isomorphic to $M_\lambda$ for some $\lambda \vdash n$ with $\ell(\lambda) \geq 3$, where $n$ is the degree of  $\overline{F(M)}$.
We expand
\begin{equation}
\overline{F(M)} = \sum_{k=0}^{\lfloor (n-1)/2 \rfloor} t_k \overline{U^n_k}.
\end{equation}
This expansion can be determined since the set $\{\overline{U^n_k} \suchthat 1 \leq k < \frac{n}{2} \}$ forms a basis for the subspace $V^n_2/\frakm^2$.
Note that $\lambda$ cannot have a pair of parts with values $k$ and $n-k$.
Using this fact together with Corollary  \ref{cor:rank2-formula-modm2},
we see that if the coefficient $t_k$ is nonnegative, then $\lambda$ has exactly $t_k$ parts with value $k$.
From this we can determine all the parts of $\lambda$ which are $< \frac{n}{2}$.
Since $\lambda$ cannot have more than one part $\geq \frac{n}{2}$, this allows us to determine the remaining part of $\lambda$, if any.
\end{proof}

\bigskip
\begin{proof}[Proof of Theorem  \ref{thm:geom}]
We write $A \sqcup B$ to denote the disjoint union of multisets $A$ and $B$.
Note that a partition may be considered to be a multiset of integers.

We  fix $n > 2$ and $\lambda \vdash n$ with $\ell(\lambda) \geq 3$,
and proceed by induction on $|J|$.
The base case $|J| = 1$ follows from Proposition \ref{prop:F-inject}.
So we assume that the statement holds for $|J| < m$ for some fixed $m > 1$.
Suppose now that
\begin{equation} \label{eqn:decomp-oFMmu}
\overline{F([M_\lambda])} = \sum_{\mu \in J} \overline{F([M_\mu])},
\end{equation}
where $|J| = m$.
Say that a pair of elements $\mu, \nu \in J$  are \emph{matching} if for some value $1 < k < n-1$ we have $k \in \mu$ and $n-k \in \nu$.
If $\mu, \nu$ are a matching pair, then we can apply Proposition \ref{prop:rank2-relations} to form a new relation of type \eqref{eqn:decomp-oFMmu}
by replacing $J$ with $J' = (J - \{\mu,\nu\}) \sqcup \{\tau\}$, where $\tau = (\mu \sqcup \nu) - \{k,n-k\}$.
At the same time, Proposition \ref{prop:rank2-relations} tells us that we also have a decomposition of base polytopes $Q(M_\tau) = Q(M_\mu) \cup Q(M_\nu)$.
Since $|J'| < m$, we can apply our induction hypothesis, and we are done.
It remains to show that there exists a matching pair in $J$.

For a partition $\tau \vdash n$, define the multiset  $g(\tau) = \{ \tau_i \suchthat \tau_i > 1,\; \tau_i  \neq \frac{n}{2} \}$. 
Define multisets $L = g(\lambda)$ and $R = \bigsqcup_{\mu  \in J} g(\mu)$.
Per Corollary  \ref{cor:rank2-formula-modm2} we expand
\[ \overline{F(M_\lambda)} =  \sum_{i=1}^{\ell(\lambda)} \overline{U^n_{\lambda_i}}, \]
and we similarly expand each $\overline{F(M_\mu)}$ on the right hand side of  \eqref{eqn:decomp-oFMmu}.
Since the set $\{\overline{U^n_k} \suchthat 1 \leq k < \frac{n}{2} \}$ forms a basis for  $V^n_2/\frakm^2$, with $\overline{U^n_k}= - \overline{U^n_{n-k}}$,
we conclude that $L \subseteq R$ and that the parts in $R-L$ can be matched into complementary pairs of the form $(k,n-k)$.
Since no partition in $J$ can contain both parts of a complementary pair, there exists a matching pair in $J$ if $R-L \neq \emptyset$.

We are assuming that $|J| \geq 2$, and that $R$ and $L$ contain all the parts not equal to $\frac{n}{2}$ or $1$ on the respective sides of \eqref{eqn:decomp-oFMmu}.
The parts equal to $1$ on both sides must match since all of the partitions have at least three parts and hence no part equal to $(n-1)$.
The only way to have $R-L = \emptyset$ 
is if there exist $\mu, \nu \in J$ each of which contains a part equal to $\frac{n}{2}$,
in which case they are matching.
Thus in all cases, there exists a matching pair $\mu, \nu \in J$,
and the result follows by induction.

\end{proof}

\noindent
Now we can give an affirmative answer to \cite[Question 7.10]{billera-jia-reiner}.

\begin{corollary} \label{cor:Hilbert-basis}
For a fixed $n$,
the Hilbert basis for the semigroup in $\Qsym/\frakm^2$ generated by the set
$S = \{ \overline{F(M_\lambda)} \suchthat \lambda \vdash n, \; \ell(\lambda) \geq 3\}$
is indexed by those $M_\lambda$ for which $\ell(\lambda) = 3$.
\end{corollary}
\begin{proof}
Let $T = \{ \overline{F(M_\lambda)} \suchthat \lambda \vdash n, \; \ell(\lambda) = 3\}$.
It follows from Proposition \ref{prop:rank2-relations} that for $\ell(\lambda) > 3$, $\overline{F(M_\lambda)}$ is decomposable into a sum $\sum_\mu \overline{F(M_\mu)}$, where for all $\mu$, $\ell(\mu) < \ell(\lambda)$. 
Hence $T$ generates the same semigroup as $S$.
As noted in \cite[Section 7]{billera-jia-reiner}, if $\ell(\lambda) = 3$, then $Q(M_\lambda)$ is indecomposable.
Theorem \ref{thm:geom} then implies that $\overline{F(M_\lambda)}$ must also be indecomposable,
so $T$ is the minimal generating set, i.e. the Hilbert basis of the semigroup.
By Proposition \ref{prop:F-inject}, distinct indexing partitions yield distinct  images,
establishing the claim.
\end{proof}

\section{Additional observations}\label{sec:observations}

In this section we discuss additional aspects of our new basis, especially regarding the expansion of $F(M)$ for a matroid $M$.

\subsection{Matroid duality, loops, and coloops}

Although we describe the basis $\{ N_\alpha \}$ as `matroid-friendly',
things are slightly less friendly when considering matroid duality in the presence of coloops.
This is due to the fact, mentioned in Section \ref{sec:loopless}, that
the mapping  $F: \Mat \to \Qsym$ factors through the quotient
$\Mat \to \Mat/\!\!\sim\; \to \Qsym$,
where $\sim$ denotes loop-coloop equivalence.


For example, a fact proved in \cite{billera-jia-reiner} is that, for any matroid $M$, in terms of the monomial basis for $\Qsym$ the following relationship holds:
\begin{equation}  \label{eqn:baseM-dual}
F(M) = \sum_\alpha m_\alpha \xx^\alpha \quad \Longrightarrow \quad F(M^{*}) = \sum_\alpha m_\alpha \xx^{\alpha^*}.
\end{equation}
where $\alpha^*$ is the \emph{reversal} of $\alpha$, obtained by writing the parts of $\alpha$ in reverse order.  
If $M$ be is a matroid of rank $r$ on $n$ elements having no loops or coloops, then we have the analogous relationship
\begin{equation}  \label{eqn:nice-dual}
F(M) = \sum_\alpha m_\alpha N_\alpha \quad \Longrightarrow \quad F(M^{*}) = \sum_\alpha m_\alpha N_{\alpha^*}.
\end{equation}
However this relationship breaks down if $M$ has loops or coloops.

We showed in Theorem  \ref{thm:FM-membership} that if $M$ is a loopless matroid of rank $r$ on $n$ elements, then $F(M) \in V^n_r$.
More generally, if $M$ is of rank $r$ on $n$ elements and has exactly $\ell$ loops, then  $F(M) \in V^n_{r+\ell}$.
Thus if $M$ has exactly $c$ coloops, then we have the duality relationship
\[ F(M) \in V^n_r \quad \Longrightarrow \quad F(M^*) \in V^n_{n-r+c}. \]

\subsection{Comultiplication}

The matroid Hopf algebra is graded by matroid rank as well as ground set size.
Let $W^n_r$ be the subspace of $\Mat$ spanned by the classes of matroids of rank $r$ on $n$ elements.
Then $W^n_r \cdot W^m_s \subset W^{n+m}_{r+s}$.
For any matroid $M$ and $A \subseteq E(M)$, $r(M) = r(M|_A) + r(M/A)$.
So comultiplication in $\Mat$ also respects these gradings. 
(For general background on Hopf algebras, see \cite{hopf-intro}.)
That is,
\begin{equation} \label{eqn:mat-coprod-grading}
\Delta W^n_r \subseteq \bigoplus_{a+b=n, \atop s+t = r} \left( W^a_s \otimes W^b_t \right).
\end{equation}

One might wonder whether the standard comultiplication of the Hopf algebra $\Qsym$ respects the grading by  the rank function for our new basis, that is, whether
\begin{equation} \label{eqn:coprod-grading}
\Delta V^n_r \subseteq \bigoplus_{a+b=n, \atop s+t = r} \left( V^a_s \otimes V^b_t \right).
\end{equation}
This is not the case.
For the simplest example, consider $n = 2$ and $r = 1$.
We have $\cN^0_0 = \{N_\fatzero\} = \{1\}$ and $\cN^2_1 = \{N_{11}\} = \{ \xx^{11} \}$.
Note that there is no $\cN^m_0$ (or rather, $\cN^m_0 = \emptyset$) for $m > 0$.
The basis vectors corresponding to the right hand side of \eqref{eqn:coprod-grading} are
\[ N_{11} \otimes N_\fatzero = \xx^{11} \otimes 1, \quad \mbox{  and  } \quad
   N_\fatzero \otimes N_{11} = 1 \otimes  \xx^{11}. 
\]
However,
\[ \Delta N_{11} = \Delta \xx^{11} = \xx^{11} \otimes 1 + \xx^1 \otimes \xx^1 + 1 \otimes  \xx^{11}, \]
which clearly does not lie in the span of the above vectors.

The failure of the comultiplication to respect the rank grading can be viewed as another artifact of loop-coloop equivalence under the morphism $F$, 
as evidenced by the fact that the rank grading \emph{is} respected by comultiplication in the quotient space corresponding to matroids with neither loops nor coloops.
Let $J \subset \Qsym$ be the ideal generated by degree one elements, i.e. by $\{N_1\} = \{\xx^1\}$.
Similarly, let $I \subset \Mat$ be the ideal generated by degree one elements, i.e. by $\{[U_{0,1}],[U_{1,1}]\}$.
Both $I$ and $J$ are Hopf ideals in their respective Hopf algebras, hence $\Mat/I$ and $\Qsym/J$ (with their naturally induced comultiplications) are Hopf algebras. 
Moreover, $I = F^{-1}(J)$, so $F: \Mat \to \Qsym$ induces a surjective Hopf algebra morphism $\Mat/I \to \Qsym/J$.
Note that a natural basis for $\Mat/I$ is the set of all matroid isomorphism classes that have neither loops nor coloops, 
while a natural basis for $\Qsym/J$ is  $\{N_\alpha \suchthat \ell(\alpha) \mbox{ is even} \}$.
Taking appropriate images under the quotient map, the relation \eqref{eqn:coprod-grading}  holds in  $\Qsym/J$.
The duality formula  \eqref{eqn:nice-dual} also holds in $\Qsym/J$.

\subsection{Comparison with other $\Qsym$ bases}

In the course of their proof in Section 10 of \cite{billera-jia-reiner}, 
the authors introduce two new $\bZ$-bases for $\Qsym$.
They also compare their bases to another $\bZ$-basis due to Stanley \cite{stanley-basis}.

Our new basis is different from these three, as evidenced by the report by those authors that all three of these bases have some negative structure constants, whereas our new basis does not.
However, of the three, ours most closely resembles that of Stanley. 
Stanley's basis element indexed by a composition $\alpha = (\alpha_1,\ldots,\alpha_m)$
is $F(P)$ where, as with our basis, $P = A_1 \oplus \cdots \oplus A_m$, is the ordered sum of antichains $A_1,\cdots,A_m$ on $\alpha_1,\ldots,\alpha_m$ elements respectively.
However, Stanley applies a natural labeling to $P$, whereas we apply an alternating labeling to the ranks in the poset for our basis.

\subsection{Surjectivity of the Hopf algebra morphism}

\BJR devote \cite[Section 10]{billera-jia-reiner} to showing that  the morphism $F:\Mat \to \Qsym$ is surjective over rational coefficients.
In this subsection we sketch one way to shorten their proof somewhat using our new basis.
The reader will need to consult \cite{billera-jia-reiner} to have the full context.

Define an ordering on compositions as follows.
To each composition $\alpha$ we assign the binary word $b(\alpha)$ that begins with $\alpha_1$ zeros followed by $\alpha_2$ ones, then $\alpha_3$ zeros, then $\alpha_4$ ones, etc.
We then linearly order compositions according to their binary words: 
$\alpha < \beta$ if $b(\alpha) <_{lex} b(\beta)$.

In their proof, \BJR make use of a novel basis for the quasisymmetric functions
based on a family of posets $\{R_\sigma\}$ of maximum rank one, indexed by binary words $\sigma \in 0\{0,1\}^{n-1}$, where $n$ is the number of elements of the poset.
We may equivalently index them using compositions of weight $n$, declaring $R_\alpha = R_{b(\alpha)}$.
We refer the reader to \cite[Section 10]{billera-jia-reiner} for the definition of this basis.
\BJR show, through a series of theorems 
that the set of $\{F(R_\alpha)\}$, where the posets are strictly labeled, forms a $\bZ$-basis for $\Qsym$.
Using our basis, one can show this more directly.
We know from  Lemma \ref{lemma:rank1} that all $\beta \in \supp(F(R_\alpha))$ are of rank $r(\alpha)$ and weight $|\alpha|$.
It is not too hard to show that the largest $\beta \in \supp(F(R_\alpha))$
with respect to the above ordering is precisely $\alpha$, and that the coefficient of $N_\alpha$ in the expansion of $F(R_\alpha)$ is one.
Thus an array giving the expansion of all the $\{F(R_\alpha)\}$ of a fixed set size $|\alpha| = n$ in terms of $\{N_\alpha\}$ with rows and columns suitably ordered is unitriangular.


\section*{Acknowledgments}

I wish to thank the referees for their helpful comments.
I extend many thanks to Isabella Novik and Sara Billey for their encouragement and their many hours of proofreading and advice through several drafts of this manuscript.
Thanks also go to Vic Reiner and Lou Billera for their clarification of parts of their paper and their helpful comments on this work.
Special thanks go to Isabella Novik for bringing the paper of Billera, Jia, and Reiner to my attention,
and to Vic Reiner for pointing me specifically to \cite[Question 7.10]{billera-jia-reiner}, which was the starting point for my investigation.

While working on this project the author was partially supported by a
graduate fellowship from VIGRE NSF Grant DMS-0354131.  
The seeds of this project were planted at the ICM in Madrid, 2006.
The trip there was made possible by funds from NSF Grant DMS-9983797.

\bibliography{kwl}
\end{document}